\documentclass[11pt]{article}
\usepackage{amsfonts,latexsym,amsmath,amscd,geometry}
\usepackage{amssymb}
\usepackage{latexsym}

\newcommand \nc{\newcommand}
\newtheorem{theorem}{Theorem}[section]
\newtheorem{lemma}[theorem]{Lemma}
\newtheorem{proposition}[theorem]{Proposition}
\newtheorem{corollary}[theorem]{Corollary}
\newtheorem{definition}[theorem]{Definition}

\newtheorem{remark}[theorem]{Remark}

\nc{\ba}{\begin{array}}\nc{\ea}{\end{array}}
\nc{\be}{\begin{eqnarray}}\nc{\ee}{\end{eqnarray}}
\nc{\beq}{\begin{equation}}\nc{\eeq}{\end{equation}}
\nc{\bex}{\begin{eqnarray*}}\nc{\eex}{\end{eqnarray*}}
\nc{\btm}{\begin{theorem}} \nc{\etm}{\end{theorem}}
\nc{\blm}{\begin{lemma}} \nc{\elm}{\end{lemma}}
\nc{\R}{\mathbb{R}} \nc{\va}{\varepsilon} \nc{\ls}{\limits}

\def\pf{\noindent{\bf Proof.\quad}}

\newcommand \qed {\hfill $\Box$}

\begin{document}

\thispagestyle{empty}
\title{Energy identity for a class of approximate biharmonic maps into sphere in dimension four}
\author{
Changyou Wang\footnote{ Department of Mathematics, University of Kentucky, Lexington, KY 40506,
cywang@ms.uky.edu.}\ \quad
Shenzhou Zheng\footnote{Department of Mathematics, Beijing Jiaotong University, Beijing 100044, P. R. China, shzhzheng@bjtu.edu.cn.}}
\date{}
\maketitle

\begin{abstract} We consider in dimension four
weakly convergent sequences of approximate biharmonic maps into sphere with bi-tension fields bounded in $L^p$ for some $p>1$.
We prove an energy identity that accounts for the loss of Hessian energies by the sum of Hessian energies over finitely many nontrivial biharmonic maps on $\mathbb R^4$.
\end{abstract}

\section {Introduction}

Let $\Omega\subset \mathbb{R}^n$ be a bounded smooth domain, and $\mathbb S^{k}\subset\mathbb R^{k+1}$ be
the unit sphere.  Recall the Sobolev space $W^{l,p}(\Omega,\mathbb S^{k})$,
$1\le l<+\infty$ and $1\le p<+\infty$, is defined by
$$
W^{l,p}(\Omega, \mathbb S^k)
=\Big\{v\in W^{l,p}(\Omega, \mathbb R^{k+1}):
\ v(x)\in \mathbb S^k \ {\rm{\ a.e.\ }}x\in\Omega\Big\}.$$

In this paper we will discuss the limiting behavior of weakly convergent sequences of approximate
(extrinsic) biharmonic maps $\displaystyle \{u_m\}\subset W^{2,2}(\Omega,\mathbb S^{k})$ in dimension $n=4$, especially an energy identity and oscillation convergence.
First we recall the notion of approximate (extrinsic) biharmonic maps.
\begin{definition} A map $u\in W^{2,2}(\Omega,\mathbb S^k)$ is called an approximate biharmonic map if there exists a bi-tension
field $h\in L^1_{\rm{loc}}(\Omega, \mathbb R^{k+1})$ such that
\begin{equation}\label{approx_biharm}
\Delta^2 u+\left(|\Delta u|^2+\Delta|\nabla u|^2+2\langle\nabla
u,\nabla\Delta u \rangle\right) u=h
\end{equation}
in the distribution sense. In particular, if $h=0$ then $u$ is called a biharmonic map to $\mathbb S^k$.
\end{definition}

Note that biharmonic maps to $\mathbb S^k$ are critical points of the Hessian energy functional
$$\displaystyle E_2(u)=\int_\Omega |\nabla^2 u|^2\,dx$$
over $W^{2,2}(\Omega, \mathbb S^k)$.
Biiharmonic maps are  higher order extensions of harmonic maps.
The study of regularity of biharmonic maps has generated considerable interests after the initial work by Chang-Wang-Yang \cite{CWY},  the readers
can refer to Wang \cite{W1,W2,W3}, Strzelecki \cite{S},  Lamm-Rivier\`e \cite{LR},
Struwe \cite{Struwe},  Scheven \cite{scheven1, scheven2} (see also Ku \cite{K} and Gong-Lamm-Wang \cite{GLW} for the boundary regularity).  In particular,  the interior regularity theorem asserts that the smoothness of $W^{2,2}$-biharmonic maps holds in dimension $n=4$,
and the partial regularity of stationary $W^{2,2}$-biharmonic maps holds in dimensions $n\ge 5$.

It is an important observation that biharmonic maps are invariant under dilations in $\mathbb R^n$ for $n=4$. Such
a property leads to non-compactness of biharmonic maps in dimension $4$, which prompts recent studies
by Wang \cite{W1} and Hornung-Moser \cite{HM} concerning the failure of strong convergence for weakly convergent biharmonic maps.  Roughly speaking,  the results in \cite{W1} and \cite{HM} assert that the failure of strong convergence occurs at finitely many concentration points of Hessian energy, where finitely many
bubbles (i.e. nontrivial biharmonic maps on $\mathbb R^4$) are generated, and the total Hessian energies
from these bubbles account for the total loss of Hessian energies during the process of convergence.

Our first result is to extend the results from \cite{W1} and \cite{HM} to the context of suitable approximate biharmonic maps to $\mathbb S^k$. More precisely,
we have
\begin{theorem}\label{energy_identity} For $n=4$, suppose $\{u_m\}\subset W^{2,2}(\Omega,\mathbb \mathbb S^{k})$ is a sequence of approximate biharmonic maps, which are bounded in $W^{2,2}(\Omega, \mathbb S^k)$ and have their bi-tension fields
$h_m$ bounded in $L^p$ for some $p>1$, i.e.
\begin{equation}\label{bound}
M:=\sup_{m}\Big(\|u_m\|_{W^{2,2}}+\|h_m\|_{L^p}\Big)<+\infty.
\end{equation}
Assume $u_m\rightharpoonup u$ in $W^{2,2}$ and $h_m\rightharpoonup h$ in $L^p$. Then \\
(i) $u$ is an approximate biharmonic map to $\mathbb S^k$ with $h$ as its bi-tension field. \\
(ii) There exist a nonnegative integer $L$ depending only on $M$
and $\{x_1,\cdots,x_L\}\subset\Omega$
such that
$\displaystyle u_m\rightarrow u \ {\rm {strongly\ in}}\
W^{2,2}_{\rm{loc}}\cap C^0_{\rm{loc}}(\Omega\backslash\{x_1,\cdots,x_L\},\mathbb S^{k}).$\\
(iii) For $1\le i\le L$, there exist a positive integer $L_i$ depending only on $M$ and nontrivial smooth biharmonic map
$\omega_{ij}$ from $\mathbb R^4$ to $\mathbb S^k$ with finite Hessian energy, $1\le j\le L_i$, such that
\begin{equation}\label{energy_id1}
\lim_{m\rightarrow\infty}\int_{B_{r_i}(x_i)} |\nabla^2u_m|^2
=\int_{B_{r_i}(x_i)}|\nabla^2 u|^2+\sum_{j=1}^{L_i}\int_{\mathbb{R}^4}|\nabla^2\omega_{ij}|^2,
\end{equation}
and
\begin{equation}\label{energy_id2}
\lim_{m\rightarrow\infty}\int_{B_{r_i}(x_i)} |\nabla u_m|^4
=\int_{B_{r_i}(x_i)} |\nabla u|^4
+\sum_{j=1}^{L_i}\int_{\mathbb{R}^4} |\nabla\omega_{ij}|^4,
\end{equation}
where $\displaystyle
r_i=\frac 12\min_{1\le j\le L,\ j\neq
i}\left\{|x_i-x_j|,\ {\rm{dist}}(x_i,\partial\Omega)\right \}.$
\end{theorem}

The idea to prove Theorem \ref{energy_identity} is based on the duality between the Lorentz spaces
$L^{2,1}$ and $L^{2,\infty}$ (see \S2 below for the definitions and basic properties). More precisely, we can bound the $L^{2,1}$-norm of $\nabla^2 u_m$ in
the neck region, while showing the $L^{2,\infty}$-norm of $\nabla^2 u_m$ can be arbitrarily small
in the neck region.  Our argument to estimate $\displaystyle\|\nabla^2 u_m\|_{L^{2,1}}$  relies heavily
on the symmetry of $\mathbb S^k$- a property that was earlier utilized by \cite{CWY}, \cite{W1}, \cite{S},
and \cite{K} in the
study of biharmonic maps. However, the argument to establish the estimation of $\displaystyle
\|\nabla^2 u_m\|_{L^{2,\infty}}$ does  not utilize the symmetry of $\mathbb S^k$ and hence holds
for any target manifold $N$.

We conjecture that Theorem 1.1 remains to be true for any target manifold $N$. For a general target manifold $N$, it remains to be a difficult question on how to obtain $L^{2,1}$-estimate for $\nabla^2 u$ similar to (2.3) in Theorem 2.3.  In a forthcoming paper \cite{WZ}, we will employ a different approach similar to \cite{HM} to prove Theorem 1.1 under the stronger assumption that the bi-tension fields are bounded in $L^p$ for some $p>\frac43$.

We would like to remark that the corresponding $L^{2,1}$-estimate
of $\nabla u$ plays a very important role in the study of harmonic maps by H\'elein \cite{H}. Later,
Lin-Rivier\`e \cite{LinR} (and Lin-Wang \cite{LW}  respectively)
utilized the duality between $L^{2,1}$ and $L^{2,\infty}$ to study the energy quantization effect of harmonic maps (and approximate harmonic maps respectively) to $\mathbb S^k$ in higher dimensions.
See also Laurain-Rivier\`e \cite{LauR} for some most recent related works.

From the view of point for applications, Theorem 1.2 can be a useful extension of the results by \cite{W1} and \cite{HM}. A typical application of Theorem 1.2 is to study asymptotic behavior at time infinity for the heat flow of biharmonic maps
in dimension $4$.

Let's review some recent studies on the heat flow of biharmonic maps undertaken by Lamm \cite{Lamm}, Gastel \cite{Gastel}, Wang \cite{W4}, and Moser \cite{M}.  For a given compact Riemannian manifold $N\subset\mathbb R^{k+1}$ without boundary,
the equation of heat flow of (extrinsic) biharmonic maps into $N$ is to seek
$u:\Omega\times [0,+\infty) \to N$ that solves (see Lamm \cite{Lamm}):
\begin{eqnarray}\label{heat_biharm}
u_t+\Delta^2 u&=&\Delta(\mathbb B(u)(\nabla u,\nabla u))+2\nabla\cdot\langle\Delta u,
\nabla(\mathbb P(u))\rangle-\langle\Delta(\mathbb P(u)),\Delta u\rangle,
\ \Omega\times (0,+\infty) \\
u&=& u_0, \ \Omega\times \{0\}\label{init_cond}\\
(u,\frac{\partial u}{\partial\nu})&=&(u_0,
\frac{\partial u_0}{\partial\nu}),
\ \partial\Omega\times (0,+\infty),\label{bdry_cond}
\end{eqnarray}
where $u_0\in W^{2,2}(\Omega, N)$ is a given map, $\mathbb P(y):\mathbb R^{k+1}\to T_y N$ is
the orthogonal projection from $\mathbb R^{k+1}$ to the tangent space of $N$ at $y\in N$,
and $\displaystyle\mathbb B(y)(X,Y)=-\nabla_X\mathbb P(y)(Y), \ \forall  X, Y\in T_y N$, is the second fundamental form of $N\subset\mathbb R^{k+1}$. Note that any $t$-independent solution
$u:\Omega\to N$ of (\ref{heat_biharm}) is a biharmonic map to $N$.

In dimension $n=4$, Lamm \cite{Lamm} established the existence of global smooth solutions to (\ref{heat_biharm})-(\ref{bdry_cond}) for $u_0\in W^{2,2}(\Omega, N)$ with small $W^{2,2}$-norm, and Gastel \cite{Gastel} and Wang \cite{W4} independently showed that
there exists a unique global weak solution to (\ref{heat_biharm}))-(\ref{bdry_cond}) for any initial data $u_0\in W^{2,2}(\Omega, N)$ that has at most finitely many singular
times. Moreover, such a solution enjoys the energy inequality:
\begin{equation}\label{energy_ineq}
2\int_0^T\int_\Omega |u_t|^2+\int_\Omega |\Delta u|^2(T)
\le \int_\Omega |\Delta u_0|^2,
\ \forall \ 0<T<+\infty.
\end{equation}
Recently, Moser \cite{M} was able to show the existence of a global weak solution to (\ref{heat_biharm})-(\ref{bdry_cond})
for any target manifold $N$ in dimensions $n\le 8$.

It follows from (\ref{energy_ineq}) that there exists a sequence $t_m\uparrow\infty$ such that $u_m:=u(\cdot, t_m)\in W^{2,2}(\Omega,N)$ satisfies \\
(i) $\displaystyle \tau_2(u_m):=\|u_t(t_m)\|_{L^2}\rightarrow 0$; and\\
(ii) $u_m$ satisfies in the distribution sense
\begin{equation}\label{approx_biharm1}
-\Delta^2 u_m+\Delta(\mathbb B(u_m)(\nabla u_m,\nabla u_m))+2\nabla\cdot\langle\Delta u_m,\nabla(\mathbb P(u_m))\rangle-\langle\Delta(\mathbb P(u_m)),\Delta u_m\rangle=\tau_2(u_m).
\end{equation}
In particular, when $N=\mathbb S^k$, by Definition 1.1 $\{u_m\}$ is a  sequence of approximate biharmonic maps to $\mathbb S^k$, which are bounded in $W^{2,2}$ and have their bi-tension fields bounded in $L^2$.

As an immediate corollary, we obtain the following theorem for the heat flow of biharmonic maps to $\mathbb S^k$  in dimension 4.
\begin{theorem}\label{energy_identity1} For $n=4$, $N=\mathbb S^k$, and $u_0\in W^{2,2}(\Omega, \mathbb S^k)$, let
$u:\Omega\times [0,+\infty)\to\mathbb S^k$, with $u\in L^\infty([0,+\infty), W^{2,2}(\Omega))$ and
$u_t\in L^2([0,+\infty), L^2(\Omega))$, be a global weak solution of (\ref{heat_biharm})-(\ref{bdry_cond}) that satisfies the
energy inequality (\ref{energy_ineq}). Then there exist $t_m\uparrow +\infty$, a biharmonic map $u_\infty\in C^\infty\cap W^{2,2}(\Omega,\mathbb S^k)$
with $u_\infty=u_0$ on $\partial\Omega$, and a nonnegative integer $L$ and  $\{x_1,\cdots, x_L\}\subset\Omega$ such that\\ (i) $u_m:=u(\cdot,t_m)\rightharpoonup u_\infty$ in $W^{2,2}(\Omega, \mathbb S^k)$.\\
(ii) $u_m\rightarrow u_\infty$ in $C^0_{\rm{loc}}\cap W^{2,2}_{\rm{loc}}(\Omega\setminus\{x_1,\cdots, x_L\}, \mathbb S^k)$.\\
(iii) For $1\le i\le L$, there exist a positive integer $L_i$ and nontrivial biharmonic maps $\{\omega_{ij}\}_{j=1}^{L_i}$
on $\mathbb R^4$ with finite Hessian energies such that
\begin{equation}\label{energy_id1}
\lim_{m\rightarrow\infty}\int_{B_{r_i}(x_i)} |\nabla^2u_m|^2
=\int_{B_{r_i}(x_i)}|\nabla^2 u_\infty|^2+\sum_{j=1}^{L_i}\int_{\mathbb{R}^4}|\nabla^2\omega_{ij}|^2,
\end{equation}
and
\begin{equation}\label{energy_id2}
\lim_{m\rightarrow\infty}\int_{B_{r_i}(x_i)} |\nabla u_m|^4
=\int_{B_{r_i}(x_i)} |\nabla u_\infty|^4
+\sum_{j=1}^{L_i}\int_{\mathbb{R}^4} |\nabla\omega_{ij}|^4,
\end{equation}
where $\displaystyle
r_i=\frac 12\min_{1\le j\le L,\ j\neq
i}\left\{|x_i-x_j|,\ {\rm{dist}}(x_i,\partial\Omega) \right\}.$
\end{theorem}

In a forthcoming article \cite{WZ}, we will show Theorem \ref{energy_identity1} remains to be true for a general target manifold $N$.

The paper is organized as follows. In \S 2, we establish the H\"older continuity for any approximate biharmonic map with its bi-tension field in $L^p$ for some $p>1$, and $L^{2,1}$-estimate of its
Hessian $\nabla^2 u$. In \S3, we show the strong convergence under  the smallness condition
of Hessian energy and set up the bubbling process. In \S4, we show no concentration of
$\displaystyle\|\nabla^2 u\|_{L^{2,\infty}(\cdot)}$  in the neck region.
In \S5, we apply the duality between $L^{2,1}$ and
$L^{2,\infty}$ to show neither Hessian energy nor oscillation  can concentrate
in the neck region, which proves Theorem \ref{energy_identity}.

\bigskip
\noindent{\bf Acknowledgement}.\  The first author is partially supported by NSF 1000115. The second author is partially supported by NSFC grant 11071012.
This work is conducted while the second author is visiting the University of Kentucky, he would like to thank Department of Mathematics for its support and hospitality.

\section{A priori estimates of approximate biharmonic maps}
\setcounter{equation}{0}
\setcounter{theorem}{0}

This section is devoted to the estimate of $L^{2,1}$ norm of
$\nabla^2 u$ for an approximate biharmonic map in terms of its Hessian energy and  $L^p$-norm of its bi-tension field for some $p>1$, and the H\"older continuity estimate under the smallness condition on its Hessian energy.

First we recall the definition and some basic properties of Lorentz spaces $L^{2,1}$ and $L^{2,\infty}$ on $\mathbb R^n$
(see \cite{H} for more detail).

\begin{definition} Let $U\subset \mathbb{R}^n$ be an open subset.
For $1<p<+\infty$ and $1\le q\le \infty$, the Lorentz space $L^{p,q}(U)$ consists of all measurable functions
$f:U\rightarrow \mathbb{R}$ such that
$$
\|f\|_{L^{p,q}(U)}=\left\{
\begin{array}{ll} \left(\int^{+\infty}_0(t^{\frac
1p}f^*(t))^q\frac{dt}t\right)^{\frac1q},\qquad \mbox{if}\ 1\le q<+\infty\\
\Big\|t^{\frac 1p}f^*(t)\Big\|_{L^{\infty}(0,+\infty)},\qquad\quad
\mbox{if}\ \ q=+\infty
\end{array}
\right.
$$
is finite, where $f^*:[0,|U|)\rightarrow \mathbb{R}$ denotes the
nonincreasing rearrangement of $|f|$:
$$
\Big|\Big\{x\in U: |f(x)|\ge s\Big\}\Big| = \Big|\Big\{t\in [0,|U|): f^*(t)\ge s\Big\}\Big|,\ \  \forall \ s\ge 0.
$$
\end{definition}

It is well-known that for $1<p<+\infty$ and $1\le q\le +\infty$,
$L^{p,q}(U)$ is the dual space of $L^{\frac{p}{p-1},\frac{q}{q-1}}(U)$. Moreover, $L^{p,p}(U)=L^p(U)$,
$L^{p',q'}(U)\subset L^{p, q}(U)$ if $1<p\le p'<+\infty$
and $1\le q'\le q\le +\infty$ and $|U|<+\infty$, and
$f\in L^{p,\infty}(U)$ is equivalent to
\begin{equation}\label{2infty}
\|f\|_{L^{p,\infty}(U)}
:=\sup_{t>0}\ t\Big|\{x\in U: |f(x)|\ge t\}\Big|^{\frac1p}<+\infty.
\end{equation}

We also recall the Sobolev embedding inequality between Lorentz spaces, whose proof can be found in \cite{H}.
\begin{proposition} \label{lorentz_embed1} For $1\le p<n$ and $1\le q<+\infty$, if
$f\in L^p(\mathbb R^n)$ and its distributional derivative
$\nabla f\in L^{p,q}(\mathbb R^n)$, then $f\in L^{\frac{np}{n-p},
q}(\mathbb R^n)$ and
\begin{equation}\label{lorentz_embed2}
\Big\|f\Big\|_{L^{\frac{np}{n-p},q}(\mathbb R^n)}
\le C\Big\|\nabla f\Big\|_{L^{p,q}(\mathbb R^n)}.
\end{equation}
\end{proposition}

Now we have
\begin{theorem}\label{21estimate} For $n=4$,
suppose  $u\in W^{2,2}(\Omega,\mathbb S^{k})$ is an approximate biharmonic map with its bi-tension field $h\in L^p(\Omega)$ for some $1<p<2$. Then for any ball $B_\delta\subset\Omega$ with radius
$\delta>0$,
$\nabla^2 u\in L^{2,1}(B_{\frac{\delta}2})$. Moreover,  for
any $0<\theta\le \frac12$ the following estimates hold
\begin{equation}\label{21bound}
\|\nabla^2 u\|_{L^{2,1}(B_{\theta\delta})}\le C\left(\theta \|\nabla^2
u\|_{L^2(B_{\delta})}+\|\nabla^2 u\|_{L^2(B_\delta)}^2+\delta^{4(1-\frac
1p)}\|h\|_{L^p(B_{\delta})}\right),
\end{equation}
and
\begin{equation}\label{22-decay}
\|\nabla^2 u\|_{L^{2}(B_{\theta\delta})}\le C\left(\theta \|\nabla^2
u\|_{L^2(B_{\delta})}+\|\nabla^2 u\|_{L^2(B_\delta)}^2+\delta^{4(1-\frac
1p)}\|h\|_{L^p(B_{\delta})}\right).
\end{equation}
\end{theorem}

\pf Since
$$\|\nabla^2 u\|_{L^2(B_{\theta\delta})}
\le \|\nabla^2 u\|_{L^{2,1}(B_{\theta\delta})},$$
(\ref{22-decay}) follows directly from (\ref{21bound}).
The idea to prove (\ref{21bound})
is similar to that of \cite{W1}. Let
$\times$ denote the wedge product in $\mathbb R^{k+1}$. First
observe that the equation (\ref{approx_biharm}) is equivalent to:
\begin{equation}\label{approx_biharm2}
\Delta\left (\nabla\cdot(\nabla u\times
u)\right)=2\nabla\cdot(\Delta u\times\nabla u)+h\times u.
\end{equation}

Since by scaling the case $\delta\not=1$ can be reduced to the case $\delta=1$,
we assume $\delta=1$ for simplicity.
Let $\widetilde{h}\in L^p(\mathbb R^4)$ be an extension of $h$ such that
\begin{equation}\label{h_extension}
\|\widetilde{h}\|_{L^q(\mathbb R^4)}
\le C\|h\|_{L^q(B_1)}, \ \forall 1\le q\le p.
\end{equation}
Let $\widetilde{u}\in W^{2,2}(\mathbb R^4,\mathbb R^{k+1})$
be an extension of $u$ such that
\begin{equation}\label{extend}
\|\nabla\widetilde u\|_{L^4(\mathbb R^n)}
\le C\|\nabla u\|_{L^{4}(B_1)},\
\|\nabla^2\widetilde u\|_{L^2(\mathbb R^n)}\le C\|\nabla^2 u\|_{L^{2}(B_1)}.
\end{equation}
Now we consider the Hodge decomposition of the 1-form $\displaystyle
d\widetilde{u}\times\widetilde{u} =\sum^4_{i=1} \frac{\partial \widetilde{u}}{\partial x_i}\times \widetilde{u} \,dx_i\in
L^4(\mathbb R^4,\wedge^1\mathbb R^4)$. It is well-known  \cite{IM}
that there exist a  function $F\in W^{1,4}(\mathbb R^4)$ and a 2-form
$H\in W^{1,4}(\mathbb R^4,\wedge^2 \mathbb R^4) $ such that
\begin{equation}\label{hodge}
d{\widetilde{u}}\times\widetilde{u} = dF+d^*H,\quad dH= 0\ \ \ \ {\rm{in}}\ \mathbb R^4,
\end{equation}
\begin{equation}\label{l4-estimate}
\left\|\nabla F\right\|_{L^4(\mathbb R^4)}+\left\|\nabla H\right\|_{L^4(\mathbb R^4)} \le
C\left\|\nabla u\right\|_{L^4(B_1)},
\end{equation}
and
\begin{equation}\label{22-estimate}
\left\|\nabla^2 F\right\|_{L^2(\mathbb R^4)}+\left\|\nabla^2 H\right\|_{L^2(\mathbb R^4)} \le
C\left\|\nabla^2 u\right\|_{L^2(B_1)},
\end{equation}
It is easy to see that $H$ satisfies
\begin{equation}\label{H-eqn}
\Delta H=d\widetilde{u}\times d\widetilde{u}\ \ {\rm{in}}\ \mathbb R^4.
\end{equation}
By Proposition \ref{lorentz_embed1} and H\"older inequality, $\displaystyle
d\widetilde{u}\times d\widetilde{u}\in L^{2,1}(\mathbb R^4)$. Hence, by the Calderon-Zgymund's $L^{p,q}$-theory, we conclude that
$\displaystyle \nabla^2 H\in L^{2,1}(\mathbb R^4)$, and
\begin{eqnarray}\label{H-21-estimate}
\|\nabla^2H\|_{L^{2,1}(\mathbb R^4)}&\le & C
\|d\widetilde{u}\times  d\widetilde{u}\|_{L^{2,1}(\mathbb R^4)}\le C \|\nabla \widetilde{u} \|^2_{L^{4,2}(\mathbb R^4)}\nonumber\\
&\le & C \|\nabla^2 \widetilde{u} \|^2_{L^{2,2}(\mathbb R^4)}
\le C \|\nabla^2 u\|^2_{L^{2}(B_1)}.
\end{eqnarray}

To estimate $\displaystyle\|\nabla^2 F\|_{L^{2,1}}$, let $G(x-y)=c_4\ln |x-y|$ be the fundamental solution of $\Delta^2$ on $\mathbb R^4$. Set
$F_1, F_2:\mathbb R^4\to \mathbb R^{k+1}$ by letting
\begin{equation}\label{F1}
F_1(x)=-2\int_{\mathbb R^4}\nabla_y G(x-y)\cdot(\Delta \widetilde{u}\times\nabla\widetilde{u})(y)\,dy,
\end{equation}
and
\begin{equation}\label{F2}
F_2(x)=\int_{\mathbb R^4} G(x-y)(\widetilde{h}\times\widetilde{u})(y)\,dy.
\end{equation}
Then it is readily seen that
\begin{equation}\label{F1_eqn}
\Delta^2 F_1=2\nabla\cdot(\Delta\widetilde{u}\times\nabla\widetilde{u}) ,
\ \ {\rm{in}}\ \ \mathbb R^4,
\end{equation}
and
\begin{equation}\label{F2_eqn}
\Delta^2 F_2=\widetilde h\times \widetilde u,
\ \ {\rm{in}}\ \ \mathbb R^4.
\end{equation}
Since $F$ satisfies
\begin{equation}\label{F_eqn}
\Delta^2 F=2\nabla\cdot(\Delta u\times \nabla u)+h\times u,
\ \ \ {\rm{in}}\ \ B_1,
\end{equation}
we conclude that $F_3:=F-F_1-F_2$ satisfies
\begin{equation}\label{F3_eqn}
\Delta^2 F_3=0, \ \ {\rm{in}}\  \ B_1.
\end{equation}

Now we want to estimate $F_1, F_2, F_3$ as follows.
For $F_1$, since
\begin{equation} \label{F1-est1}
\nabla^3 F_1(x)=
 -2\int_{\mathbb R^4} \nabla^4_y G(x-y) \cdot(\Delta\widetilde u\times \nabla\widetilde u)(y)\,dy,
\end{equation}
and by H\"older inequality $\Delta\widetilde u\times \nabla\widetilde u\in L^{\frac43,1}(\mathbb R^4)$,
we have by Calderon-Zygmund $L^{p,q}$-theory that $\nabla^3 F_1\in L^{\frac43,1}(\mathbb R^4)$ and
\begin{equation}\label{F1-est2}
\Big\|\nabla^3 F_1\Big\|_{L^{\frac43,1}(\mathbb R^4)}
\le C\Big\|\Delta\widetilde u\times \nabla\widetilde u\Big\|_{L^{\frac43,1}(\mathbb R^4)}
\le C\Big\|\nabla^2 \widetilde u\Big\|_{L^2(\mathbb R^4)}
\Big\|\nabla\widetilde u\Big\|_{L^{4,2}(\mathbb R^4)}
\le C\Big\|\nabla^2 u\Big\|_{L^2(B_1)}^2.
\end{equation}
Hence, by Proposition 2.2 we have that $\nabla^2 F_1\in L^{2,1}(\mathbb R^4)$ and
\begin{equation}\label{F1-est3}
\Big\|\nabla^2 F_1\Big\|_{L^{2,1}(\mathbb R^4)}
\le C\Big\|\nabla^3 F_1\Big\|_{L^{\frac43,1}(\mathbb R^4)}
\le C\Big\|\nabla^2 u\Big\|_{L^2(B_1)}^2.
\end{equation}
For $F_2$, we have
\begin{eqnarray*}
\Big|\nabla^2 F_2(x)\Big|&\le& c_4\Big|\int_{\mathbb R^4} |x-y|^{-2}(\widetilde h\times\widetilde u)(y)\,dy\Big|\le  CI_2(|\widetilde h|)(x)
\end{eqnarray*}
where $I_\beta(f)$ is the Riesz potential of order $\beta$ ($0<\beta\le 4$) defined by
\begin{equation} \label{reisz}
I_{\beta}(f)(x)\equiv\int_{\mathbb{R}^4}\frac
{f(y)}{|x-y|^{4-\beta}} \,dy, \qquad x\in \mathbb{R}^4.
\end{equation}
It follows from Adams \cite{A} (see also \cite{W1})
that $\nabla^2 F_2\in
L^{\frac{2p}{2-p}}(\mathbb R^4)$ and
\begin{equation}\label{F2-est1}
\Big\|\nabla^2 F_2\Big\|_{L^{\frac{2p}{2-p}}(\mathbb R^4)}
\le C\Big\|I_2(\widetilde h)\Big\|_{L^{\frac{2p}{2-p}}(\mathbb R^4)}\le C\Big\|\widetilde h\Big\|_{L^p(\mathbb R^4)}
\le C\Big\|h\Big\|_{L^p(B_1)}.
\end{equation}
Since $L^{\frac{2p}{2-p}}(B_1)\subset L^{2,1}(B_1)$,
this implies
\begin{equation}\label{F2-est2}
\Big\|\nabla^2 F_2\Big\|_{L^{2,1}(B_1)}
\le C\Big\|\nabla^2 F_2\Big\|_{L^{\frac{2p}{2-p}}(B_1)}
\le C\Big\|h\Big\|_{L^p(B_1)}.
\end{equation}
Since $F_3$ is a biharmonic function on $B_1$, the standard estimate implies
that for any $0<\theta\le\frac12$,
\begin{eqnarray}\label{F3-est1}
\Big\|\nabla^2 F_3\Big\|_{L^{2,1}(B_\theta)}
&\le& C\theta \Big\|\nabla^2 F_3\Big\|_{L^2(B_1)}\nonumber\\
&\le& C\theta\Big(\|\nabla^2 F\|_{L^2(B_1)}
+\|\nabla^2 F_1\|_{L^2(B_1)}
+\|\nabla^2 F_2\|_{L^2(B_1)}\Big)\nonumber\\
&\le& C\theta \Big(\|\nabla^2 u\|_{L^2(B_1)}+\|h\|_{L^p(B_1)}\Big).
\end{eqnarray}
Combining (\ref{F1-est3}), (\ref{F2-est1}) together with  (\ref{F3-est1})
yields
\begin{equation}\label{F-est}
\Big\|\nabla^2 F\Big\|_{L^{2,1}(B_\theta)}
\le C\Big(\theta\Big\|\nabla^2 u\Big\|_{L^2(B_1)}+
\Big\|h\Big\|_{L^p(B_1)}\Big).
\end{equation}
Since $|u|=1$, we have
\begin{equation}\label{point_id}
\frac{\partial^2 u}{\partial x_i\partial x_j}\cdot u+\frac{\partial u}{\partial x_i}\cdot\frac{\partial u}{\partial x_j}=0,
\end{equation}
and hence
\begin{eqnarray}\label{21compare}
\left|\frac{\partial^2 u}{\partial x_i\partial x_j}\right|
&\le&\left|\frac{\partial^2 u}{\partial x_i\partial x_j}\cdot u\right|
+\left|\frac{\partial^2 u}{\partial x_i\partial x_j}\times u\right|\nonumber\\
&\le& 2\left|\frac{\partial u}{\partial x_i}\right|
\left|\frac{\partial u}{\partial x_j}\right|+\left|\frac{\partial}{\partial x_i}
(\frac{\partial u}{\partial x_j}\times u)\right|\\
&\le& 2|\nabla u|^2 +C(|\nabla^2 F|+|\nabla^2 H|),
\end{eqnarray}
and by (\ref{H-21-estimate}) and (\ref{F-est}) we  have
\begin{eqnarray*}
\Big\|\nabla^2 u\Big\|_{L^{2,1}(B_{\theta})}
&\le& C\Big(\Big\||\nabla u|^2\Big\|_{L^{2,1}(B_{\theta})}
+\Big\|\nabla^2 F\Big\|_{L^{2,1}(B_{\theta})}
+\Big\|\nabla^2 H\Big\|_{L^{2,1}(B_{\theta})}\Big)\\
&\le&
C\Big(\theta \Big\|\nabla^2 u\Big\|_{L^2(B_1)}
+\Big\|\nabla^2 u\Big\|_{L^2(B_1)}^2+\Big\|h\Big\|_{L^p(B_1)}\Big).
\end{eqnarray*}
This clearly yields (\ref{21bound}). \qed\\

As a corollary of Theorem \ref{21estimate} and (\ref{22-decay}), we can show that an approximate biharmonic map to $\mathbb S^k$, with the bi-tension field in $L^p$ for some $1<p<2$, is H\"older
continuous. More precisely, we have
\begin{corollary} \label{holder}
For $n=4$, there exist an $\epsilon_0>0$,
$\alpha\in (0,1)$, and $2<q\le \frac{2p}{2-p}$ such that if $u\in W^{2,2}(\Omega,\mathbb S^k)$ is an approximate biharmonic map with its bi-tension field $h\in L^p(\Omega)$ for
some $1<p<2$, that satisfies
\begin{equation}\label{smallness1}
\Big\|\nabla^2 u\Big\|_{L^2(B_r(x))}\le\epsilon_0,
\end{equation}
for some $B_r(x)\subset\Omega$, then \\
(i) $u\in C^\alpha(B_r(x),\mathbb S^k)$ and
\begin{equation}\label{holder_est}
\Big[u\Big]_{C^\alpha(B_{\frac{r}2}(x))}
\le C\Big(\epsilon_0, \|h\|_{L^p(B_r(x))}\Big),
\end{equation}
(ii) $\nabla^2 u\in L^q(B_\frac{r}2)$ and
\begin{equation}\label{q_bound}
\Big\|\nabla^2 u\Big\|_{L^q(B_\frac{r}2)}\le
C\Big(\epsilon_0, \|h\|_{L^p(B_r(x))}\Big).
\end{equation}
\end{corollary}
\pf For simplicity, assume $x=0$ and $r=1$. For any ball
$B_\delta(y)\subset B_1\subset\Omega$, (\ref{smallness1})
yields
$$\Big\|\nabla^2 u\Big\|_{L^2(B_\delta(y))}\le\epsilon_0$$
so that applying (\ref{22-decay}),
we have that for any $\theta\in (0,\frac12)$,
\begin{equation}\label{22-decay1}
\Big\|\nabla^2 u\Big\|_{L^2(B_{\theta\delta}(y))}
\le C(\theta+\epsilon_0)\Big\|\nabla^2 u\Big\|_{L^2(B_\delta(y))}
+C\delta^{4(1-\frac{1}{p})}\Big\|h\Big\|_{L^p(B_\delta(y))}.
\end{equation}
Thus by choosing sufficiently small $\epsilon_0>0$ and $\theta=\theta_0\in (0,\frac12)$ we have
\begin{equation}\label{22-decay2}
\Big\|\nabla^2 u\Big\|_{L^2(B_{\theta_0\delta}(y))}
\le \frac12\Big\|\nabla^2 u\Big\|_{L^2(B_\delta(y))}
+C\delta^{4(1-\frac{1}{p})}\Big\|h\Big\|_{L^p(B_\delta(y))}.
\end{equation}

It is standard that iterations of (\ref{22-decay1}) would imply that there exists $\alpha\in (0,1)$ such that for any $y\in
B_\frac12$ and $0<R\le\frac14$,
\begin{equation}
\int_{B_R(y)}|\nabla^2 u|^2
\le CR^{2\alpha}\Big(\int_{B_1}|\nabla^2 u|^2+
\|h\|_{L^p(B_1)}^2\Big).
\end{equation}
This, combined with  Morrey's decay Lemma, clearly implies
that $u\in C^\alpha(B_\frac12)$ and
$$\Big[u\Big]_{C^\alpha(B_\frac12)}
\le C\Big(\Big\|\nabla^2 u\Big\|_{L^2(B_1)}+
\Big\|h\Big\|_{L^p(B_1)}\Big).$$
This proves (i).

(ii) can be proved by applying  Adams's Riesz potential estimate
between Morrey spaces. First recall that for an open set
$U\subset\mathbb R^4$, $1\le p<+\infty$, $0<\lambda\le 4$, the
Morrey space $M^{p,\lambda}(U)$ is defined by
$$
M^{p,\lambda}(U)
=\Big\{f\in L^p(U):
\|f\|_{M^{p,\lambda}}^p
=\sup_{B_r\subset U} r^{\lambda-4}\int_{B_r}|f|^p<+\infty \Big\}.
$$
By (i), we know that for some $0<\alpha<1$,
$\nabla^2 u\in M^{2,4-2\alpha}(B_\frac34)$ and
\begin{equation}\label{morrey1}
\Big\|\nabla^2 u\Big\|_{M^{2,4-2\alpha}(B_\frac34)}
\le C\Big(\|\nabla^2 u\|_{L^2(B_1)}+\|h\|_{L^p(B_1)}\Big).
\end{equation}
Let $\widetilde u:\mathbb R^4\to \mathbb R^{k+1}$ be an extension
of $u$ on $B_\frac34$ such that
\begin{equation}\label{u_extension}
\Big\|\nabla^2 \widetilde u\Big\|_{M^{2,4-2\alpha}(\mathbb R^4)}
\le C\Big\|\nabla^2 u\Big\|_{M^{2,4-2\alpha}(B_\frac34)}.
\end{equation}
Then there exists a harmonic function $v: B_\frac34\to\mathbb R^{k+1}$,  with $\|\nabla^2 v\|_{L^2(B_\frac34)}\le C\|\nabla^2 u\|_{L^2(B_1)}$, such that $u=v+w$ on $B_\frac34$, where
$$
w(x)=\int_{\mathbb R^4} \frac{\Delta \widetilde u(y)}{|x-y|^2} \,dy.$$
Hence
$$|\nabla w|(x)\le CI_1(|\Delta \widetilde u|)(x), \ x\in\mathbb R^4.$$
By Adams's Riesz potential estimate (see \cite{A} or \cite{W1}), we have $|\nabla w|\in M^{\frac{4-2\alpha}{1-\alpha},
4-2\alpha}(\mathbb R^4)$.
Therefore we can conclude that
$\nabla u\in M^{\frac{4-2\alpha}{1-\alpha}, 4-2\alpha}(B_\frac58)$
and
\begin{equation} \label{gradient_est}
\Big\|\nabla u\Big\|_{L^{\frac{4-2\alpha}{1-\alpha}}(B_\frac58)}
\le \Big\|\nabla u\Big\|_{M^{\frac{4-2\alpha}{1-\alpha},4-2\alpha}(\mathbb R^4)}
\le C\Big(\|\nabla^2 u\|_{L^2(B_1)}+\|h\|_{L^p(B_1)}\Big).
\end{equation}
Denote $q_1=\frac{4-2\alpha}{1-\alpha}>4$. Then we have
$\Delta u\times \nabla u\in L^{\frac{2q_1}{2+q_1}}(B_\frac58)$. Note
that $q_2=\frac{2q_1}{2+q_1}>\frac43$.  Let $C:B_\frac58\to\mathbb R^{k+1}$ be the
harmonic function extension  of $\nabla\cdot(\nabla u\times u)$ from $\partial B_\frac58$
to $B_\frac58$.
Then,  by the equation (\ref{approx_biharm2}), we have that  in
$B_\frac58$,
\begin{eqnarray*}
\nabla\cdot(\nabla u\times u)(x)
&=&-2\int_{B_\frac58}\nabla_y \hat{G}(x-y)(\Delta u\times \nabla u)(y)\,dy
+\int_{B_\frac58} \hat{G}(x-y) (h\times u)(y)\,dy+C(x)\\
&=& A(x)+B(x)+C(x),
\end{eqnarray*}
where $\hat{G}$ is the Green function of $\Delta$ on $B_\frac58$.
It is easy to see that
$$|A(x)|\le C|I_1(x)|
\le I_1\Big(|\Delta u||\nabla u|\chi_{B_\frac58}\Big)(x), \ x\in B_\frac58,$$
where $\chi_{B_\frac58}$ is the characteristic function of $B_\frac58$. Hence we have
$A(x)\in L^{\frac{4q_2}{4-q_2}}(B_\frac58)$, and
$$\Big\|A(x)\Big\|_{L^{\frac{4q_2}{4-q_2}}(B_\frac58)}
\le C\Big\||\Delta u||\nabla u|\Big\|_{L^{q_2}(B_\frac58)}
\le C\Big(\|\nabla^2 u\|_{L^2(B_1)}+\|h\|_{L^p(B_1)}\Big).
$$
For $B(x)$, it is easy to see that
$$|B(x)|\le CI_2(|h|)(x),
\ \forall x\in B_\frac58.$$
Hence $B(x)\in L^{\frac{2p}{2-p}}(B_\frac58)$, and
$$\Big\|B(x)\Big\|_{L^{\frac{2p}{2-p}}(B_\frac58)}
\le C\Big\|h\Big\|_{L^p(B_1)}.$$
It is easy to see that
$$\Big\|C(x)\Big\|_{L^\infty(B_\frac12)}
\le\Big\|\nabla\cdot(\nabla u\times u)\Big\|_{L^2(\partial B_\frac58)}
\le C\Big\|\nabla^2 u\Big\|_{L^2(B_1)}.$$
Set $q=\min\{\frac{2p}{2-p}, \ \frac{4q_2}{4-q_2}\}$. Since $q_2>\frac43$, we have
$q>2$ and
$$\Big\|\nabla\cdot(\nabla u\times u)\Big\|_{L^{q}(B_\frac12)}
\le C\Big(\|\nabla^2 u\|_{L^2(B_1)}+\|h\|_{L^p(B_1)}\Big).$$
This, combined with (\ref{21compare}) and (\ref{gradient_est}), implies that
(\ref{q_bound}) holds for some $q>2$. The proof is now complete.
\qed

\section {Blow up analysis and energy inequality}
\setcounter{equation}{0}
\setcounter{theorem}{0}

This section is devoted to $\epsilon_0$-compactness lemma and
preliminary steps on the blow up analysis of approximate
biharmonic maps with bi-tension fields bounded in $L^p$ for
some $p>1$. In particular, we will indicate
that (\ref{energy_id1}) holds with ``$=$ '' replaced by ``$\ge$''.

First we have
\begin{lemma}\label{e-strong} For $n=4$,
there exists an $\epsilon_0>0$ such that if $\{u_m\}\subset
W^{2,2}(B_1,\mathbb S^k)$ is a sequence of approximate biharmonic
maps satisfying
\begin{equation}\label{smallness2}
\sup_{m}\left\|\nabla^2 u_m\right\|_{L^2(B_1)}\le\epsilon_0,
\end{equation}
and $u_m\rightharpoonup u$ in $W^{2,2}(B_1)$ and
$h_m\rightharpoonup h$ in $L^p(B_1)$ for some $p>1$.
Then $u$ is an approximate biharmonic map with bi-tension field
$h$, and
\begin{equation}\label{strong}
\lim_{m\rightarrow\infty}
\Big\|u_m-u\Big\|_{W^{2,2}(B_\frac12)}=0.
\end{equation}
\end{lemma}
\pf The first assertion follows easily from (\ref{approx_biharm2}). To show (\ref{strong}), it suffices
to show that $\{u_m\}$ is a Cauchy sequence in $W^{2,2}(B_\frac12)$.
By (\ref{smallness2}) and Corollary \ref{holder}, there exist
$\alpha\in (0,1)$ and $q>2$ such that
$$\sup_m\Big[ \Big\|u_m\Big\|_{C^\alpha(B_\frac34)}
+\Big\|\nabla^2 u_m\Big\|_{L^q(B_\frac34)}\Big]\le C.$$
Hence we may assume that
$$\lim_{m, l\rightarrow\infty} \Big\|u_m-u_l\Big\|_{L^\infty(B_\frac34)}=0.$$
For $\eta\in C_0^\infty(B_\frac34)$ be a cut-off function of
$B_\frac12$, multiplying the equations of $u_m$ and $u_l$ by
$(u_m-u_l)\phi^2 $ and integrating over $B_1$, we obtain
\begin{eqnarray*}
&&\int_{B_1}|\Delta(u_m-u_l)|^2\phi^2\\
&\le &\int_{B_1}|\Delta(u_m-u_l)|(2|\nabla(u_m-u_l)||\nabla\phi^2|+|u_m-u_l||\Delta\phi^2|)
+\int_{B_1}|h_m-h_l||u_m-u_l|\phi^2\\
&&+3\int_{B_1}(|\Delta u_l|^2 +|\Delta u_m|^2)|u_m-u_l|\phi^2\\
&&+4\int_{B_1}|\nabla^2 u_m||\nabla u_m||\nabla(u_m(u_m-u_l)\phi^2)|\\
&&+4\int_{B_1}|\nabla^2 u_l||\nabla u_l||\nabla(u_l(u_m-u_l)\phi^2)|\\
&=&I+II+III+IV+V.
\end{eqnarray*}
It is easy to see
$$|I|\le C(\|\nabla (u_m-u_l)\|_{L^2(B_\frac34)}
+\|u_m-u_l\|_{L^\infty(B_\frac34)})\rightarrow 0,$$
$$|II|\le C\|h_m-h_l\|_{L^1(B_\frac 34)}\|u_m-u_l\|_{L^\infty(B_\frac34)}\rightarrow 0,$$
$$|III|\le C(\|\nabla^2 u_m\|_{L^2(B_\frac34)}^2
+\|\nabla^2 u_l\|_{L^2(B_\frac34)}^2)\|u_m-u_l\|_{L^\infty(B_\frac34)}\rightarrow 0.$$
For $IV$, observe that for $1<r<4$ with $\frac14+\frac1{q}
+\frac1{r}=1$, we have
\begin{eqnarray*}
|IV|
&\le& C\Big(\|\nabla^2 u_m\|_{L^2(B_\frac34)}\|\nabla u\|_{L^4(B_\frac34)}^2\|u_m-u_l\|_{L^\infty(B_\frac34)}\\
&&+\|\nabla^2 u_m\|_{L^q(B_\frac34)}\|\nabla u_m\|_{L^4(B_\frac34)}
\|\nabla(u_m-u_l)\|_{L^r(B_\frac34)}\Big)
\rightarrow 0,
\end{eqnarray*}
since $\displaystyle\|\nabla(u_m-u_l)\|_{L^r(B_\frac34)}\rightarrow 0$.
Similarly, we can show
$$|V|\rightarrow 0.$$
Hence $\{u_m\}$ is a Cauchy sequence in $W^{2,2}(B_\frac12)$.
This completes the proof.
\qed
\begin{lemma} Under the same assumptions as Theorem \ref{energy_identity}, there exists a finite subset $\Sigma\subset
\Omega$ such that
$u_m\rightarrow u$ in $W^{2,2}_{\rm{loc}}\cap C^0_{\rm{loc}}
(\Omega\setminus\Sigma,\mathbb S^k)$. Moreover,
$u\in W^{2,2}\cap C^0(\Omega,\mathbb S^k)$ is an approximate
biharmonic map with bi-tension field $h$.
\end{lemma}
\pf Let $\epsilon_0>0$ be given by Corollary 2.4, and define
\begin{equation}\label{concentration}
\Sigma:=\bigcap_{r>0}\Big\{x\in\Omega:
\liminf_{m\rightarrow\infty}\int_{B_r(x)}|\nabla^2 u_m|^2>\epsilon_0^2\Big\}.
\end{equation}
Then by a simple covering argument we have that
$\Sigma$ is a finite set. In fact
$$H^0(\Sigma)\le\frac{1}{\epsilon_0^2}\sup_{m}\int_\Omega
|\nabla^2 u_m|^2<+\infty.$$
For any $x_0\in\Omega\setminus\Sigma$, there exists $r_0>0$ such
that
$$\liminf_{m\rightarrow\infty}\int_{B_{r_0}(x_0)}
|\nabla^2 u_m|^2\le\epsilon_0^2.$$
Hence Corollary 2.4 and Lemma 3.1 imply that there exists
$\alpha\in (0,1)$ such that
$$\Big\|u_m\Big\|_{C^\alpha(B_{\frac{r_0}2}(x_0))}
\le C,$$
so that $u_m\rightarrow u$ in $C^0\cap W^{2,2}(B_{\frac{r_0}2}(x_0))$. This proves that $u_m\rightarrow u$ in
$W^{2,2}_{\rm{loc}}\cap C^0_{\rm{loc}}(\Omega\setminus\Sigma)$.
It is clear that $u\in W^{2,2}(\Omega)$ is an approximate biharmonic map with bi-tension field $h\in L^p(\Omega)$. Applying Corollary 2.4 again, we conclude that $u\in C^0(\Omega,\mathbb S^k)$. \qed\\

\noindent{\bf Proof of Theorem}\ref{energy_identity}:

\smallskip

The proof of (\ref{energy_id1}) with ``$\ge$'' is similar to
\cite{W1} Lemma 3.3. Here we sketch it. For any $x_0\in\Sigma$,
there exist $r_0>0$, $x_m\rightarrow x_0$ and $r_m\downarrow 0$ such that
$$
\max_{x\in B_{r_0}(x_0)}
\Big\{\int_{B_{r_m}(x)}|\nabla^2 u_m|^2
\Big\}
=\frac{\epsilon_0^2}2=\int_{B_{r_m}(x_m)}|\nabla^2 u_m|^2.
$$
Define $v_m(x)=u_m(x_m+r_mx): B_{\frac{r_0}{r_m}}\to \mathbb S^k$. Then $v_m$ is an approximate biharmonic map, with bi-tension
field $\widetilde{h_m}(x)=r_m^4 h(x_m+r_m x)$, that satisfies
$$\int_{B_1(x)}|\nabla^2 u_m|^2\le\frac{\epsilon_0^2}2,
\ \forall x\in B_{\frac{r_0}{r_m}},
\ {\rm{and}}\ \int_{B_1}|\nabla^2 u_m|^2=\frac{\epsilon_0^2}2,$$
and
$$\Big\|\widetilde{h_m}\Big\|_{L^p(B_{\frac{r_0}{r_m}})}
\le r_m^{4(1-\frac1p)}\Big\|h_m\Big\|_{L^p(\Omega)}
\rightarrow 0.$$
Thus Corollary 2.4 and Lemma 3.1 imply that,
after taking possible subsequences, there exists a nontrivial biharmonic map $\omega:\mathbb R^4\to\mathbb S^k$ with
$$\frac{\epsilon_0^2}2\le \int_{\mathbb R^4}|\nabla^2 \omega|^2<+\infty
$$
such that $v_m\rightarrow \omega$ in
$W^{2,2}_{\rm{loc}}\cap C^0_{\rm{loc}}(\mathbb R^4)$.
Performing such a blow-up argument near any $x_i\in\Sigma$,
$1\le i\le L$, we can find all possible nontrivial biharmonic maps
$\{\omega_{ij}\}\in W^{2,2}(\mathbb R^4)$ for $1\le j\le L_i$,
with $L_i\le CM\epsilon_0^{-2}$. It is not hard to see
(\ref{energy_id1}) holds with ``$=$'' replaced by ``$\ge$''.

In order prove ``$\le$''  of (\ref{energy_id1}),
we need to show that the $L^{2,\infty}$-norm of $u_m$ over any
neck region is arbitrarily small. This will be done in the next
section. We will return to the proof of Theorem \ref{energy_identity} in \S5.

\section{$L^{2,\infty}$-estimate in the neck region}
\setcounter{equation}{0}
\setcounter{theorem}{0}

In this section, we first show that there is no concentration of
$\|\nabla^2 u\|_{L^{2,\infty}}$ in the neck region. Then use
the duality between $L^{2,1}$ and $L^{2,\infty}$ to prove Theorem \ref{energy_identity} by showing that there is no Hessian energy concentration in the neck region.
More precisely, we have

\begin{lemma}\label{2infty} For any $\epsilon>0$, suppose that $u\in W^{2,2}(B_1,\mathbb S^k)$
is an approximate biharmonic map,
whose bi-tension field $h\in L^p(B_1)$ for some $p>1$,
satisfying
that for $0<\delta<\frac12$, $R>1$, and $0<r<\frac{2\delta}{R}$,
\begin{equation}\label{small_energy}
\int_{B_{2\rho}\backslash B_\rho}|\nabla^2u|^2\,dx \le
\epsilon^2,\qquad \forall \  Rr\le \rho \le 2\delta,
\end{equation}
then
\begin{eqnarray}\label{small_2infty}
\Big\|\nabla^2 u\Big\|_{L^{(2,\infty)}(B_{{\delta}\backslash B_{2Rr}})} &\le&
C\Big[\epsilon+\epsilon^{\frac12}\|\nabla^2 u\|_{L^2(B_{2\delta})}
+\delta^2+\delta^{4(1-\frac 1p)}\Big\|h\Big\|_{L^p(B_{2\delta})}
\nonumber\\
&&+\epsilon^{-1}\delta^{8(1-\frac1p)}\|h\|_{L^p(B_{2\delta})}^2\Big].
\end{eqnarray}
\end{lemma}
\pf First recall that
\begin{eqnarray}
\Big\|\nabla^2 u\Big\|_{L^{2,\infty}(B_{\delta}\backslash B_{2Rr})}^2
&=&
\sup_{\lambda>0}\ \lambda^2\Big|\Big\{x\in B_{\delta}\backslash
B_{2Rr}:\ |\nabla^2 u|>\lambda\Big\}\Big|\nonumber\\
&\le &
\sup_{0<\lambda\le 1}\ \lambda^2\Big|\Big\{x\in
B_{\delta}\backslash B_{2Rr}:\ |\nabla^2
u|>\lambda\Big\}\Big|\nonumber\\
&+&
\sup_{\lambda>1}\ \lambda^2\Big|\Big\{x\in B_{\delta}\backslash
B_{2Rr}:\ |\nabla^2 u|>\lambda\Big\}\Big|\nonumber\\
&\le & C \delta^4+\sup_{\lambda>1}\  \lambda^2\Big|\Big\{x\in
B_{\delta}\backslash B_{2Rr}:\ |\nabla^2 u|>\lambda\Big\}\Big|.
\end{eqnarray}
It suffices to estimate
$\displaystyle \lambda^2\Big|\Big\{x\in
B_{\delta}\backslash B_{2Rr}:\ |\nabla^2 u|>\lambda\Big\}\Big|$
for $\lambda>1$.

We may assume that
$\delta =2^{K} \sqrt{\frac{\epsilon}{\lambda}}$ for some positive
integer $K\ge 1$. There are two cases to consider:\\
(i) $\sqrt{\frac{\epsilon}{\lambda}}\ge 2Rr$. Then we have
\begin{eqnarray*}
&&\lambda^2\Big|\Big\{x\in
B_{\delta}\backslash B_{2Rr}:\ |\nabla^2 u|>\lambda\Big\}\Big|\\
&\le&  \lambda^2\Big|\Big\{x\in
B_{\sqrt{\frac{\epsilon}{\lambda}}}\backslash B_{2Rr}:\
|\nabla^2 u|>\lambda\Big\}\Big|
+\lambda^2\Big|\Big\{x\in
B_{\delta}\backslash B_{\sqrt{\frac{\epsilon}{\lambda}}}:\
|\nabla^2 u|>\lambda\Big\}\Big|\\
&\le& C\epsilon^2
+\lambda^2\sum_{i=0}^{K-1}\Big|\Big\{x\in
B_{2^{i+1}\sqrt{\frac{\epsilon}{\lambda}}}\backslash B_{2^i\sqrt{\frac{\epsilon}{\lambda}}}
: |\nabla^2 u|>\lambda\Big\}\Big|.
\end{eqnarray*}
(ii) $\sqrt{\frac{\epsilon}{\lambda}}<2Rr$. Then we may assume that
there exists $1\le i_0\le K$ such that
$2Rr=2^{i_0}\sqrt{\frac{\epsilon}{\lambda}}$ so that
\begin{eqnarray*}
\lambda^2\Big|\Big\{x\in
B_{\delta}\backslash B_{2Rr}:\ |\nabla^2 u|>\lambda\Big\}\Big|
&\le&  \lambda^2\sum_{i=i_0}^{K-1}
\Big|\Big\{x\in B_{2^{i+1}\sqrt{\frac{\epsilon}{\lambda}}}
\setminus
B_{2^i\sqrt{\frac{\epsilon}{\lambda}}}:\ |\nabla^2 u|>\lambda\Big\}\Big|.
\end{eqnarray*}

It is not hard to see that the case (ii) can be done by the same
way as the case (i). Thus we only need to prove (i). To simplify
the presentation, introduce
$$r_0:=\sqrt{\frac{\epsilon}{\lambda}},
\ {\mathcal A_i}:=B_{2^{i+1}r_0}\setminus B_{2^i r_0}
\ {\rm{and}}\ {\mathcal B_i}:=B_{2^{i+2}r_0}\setminus B_{2^{i-1} r_0},
\ 0\le i\le K-1.$$
For $0\le i\le K-1$,  let $u_i: \mathbb R^4\to \mathbb R^{k+1}$ be an extension of $u$ from
${\mathcal B_i}$ to $\mathbb R^4$ such that
\begin{equation} \label{u_extension}
|u_i|\le 1;\  \int_{\mathbb{R}^4}(|\nabla^2 u_i|^2
+|\nabla u_i|^4)\,dx\le C
\int_{\mathcal B_i} (|\nabla^2 u|^2 +|\nabla u|^4)\,dx\le C\int_{\mathcal B_i}|\nabla^2 u|^2\,dx
\le C \epsilon^2,
\end{equation}
where we have used $|\nabla u|^2=|\Delta u\cdot u|\le |\Delta u|$ on $\mathcal B_i$ and (\ref{small_energy}).
Let $h_i:\mathbb R^4\to\mathbb R^{k+1}$ be an extension of $h$ from $\mathcal B_i$ to $\mathbb R^4$ such that
\begin{equation}\label{h_extension}
\Big\|h_i\Big\|_{L^q(\mathbb R^4)}
\le C\Big\| h\Big\|_{L^q(\mathcal B_i)},
\ 1\le q\le p.
\end{equation}
Define $v_i:\mathbb R^4\to\mathbb R^{k+1}$ by
\begin{equation}\label{v_i}
v_i(x) = \int_{\mathbb{R}^4} G(x-y)\left[h_i- (|\Delta u_i|^2+\Delta(|\nabla u_i|^2)+2\langle\nabla u_i,\nabla\Delta u_i\rangle
)u_i\right](y)\, dy,
\end{equation}
where $G(\cdot)$ is the fundamental solution of $\Delta^2$ on $\mathbb R^4$,
and $w_i: \mathcal B_i\to\mathbb R^{k+1}$ by
\begin{equation} \label{w_i}
w_i=u-v_i.
\end{equation}
Then it is easy to see that
\begin{equation}\label{eqn_wi}
\Delta^2w_i=0\
\ \ {\rm{on}}\ \ \mathcal B_i.
\end{equation}

Now we have\\
\noindent {\bf Claim}. For $0\le i\le K-1$,
\begin{equation}\label{vi-estimate}
\Big\|\nabla^2 v_i\Big\|_{L^{2,\infty}(\mathbb R^4)}
\le C\Big[
\Big\|h\Big\|_{L^1(\mathcal B_i)}+\Big\|\nabla^2 u\Big\|_{L^2(\mathcal B_i)}^{\frac32}
\Big].
\end{equation}
and
\begin{equation}\label{wi-estimate}
\sup_{x\in \mathcal A_i}\Big\{|x|^2\Big|\nabla^2 w_i\Big|(x)\Big\}
\le  C\Big[ \Big\|h\Big\|_{L^1(\mathcal B_i)} +\Big\|\nabla^2 u\Big\|_{L^2(\mathcal B_i)}\Big].
\end{equation}
To see (\ref{vi-estimate}), note that
\begin{eqnarray}\label{vi-estimate1}
\nabla^2 v_i(x)&=&\int_{\mathbb{R}^4} \nabla^2_y G(x-y)
\left[h_i-
(|\Delta u_i|^2+\Delta(|\nabla u_i|^2)+2\langle\nabla u_i,
\nabla\Delta u_i\rangle )u_i\right](y)\,dy\nonumber\\
&=&I(x)+II(x)+III(x)+IV(x).
\end{eqnarray}
Since $|\nabla^2_y G(x-y)|\le \frac{C}{|x-y|^2}\in
L^{2,\infty}(\mathbb{R}^4)$, we have
$$\Big\|I\Big\|_{L^{2,\infty}(\mathbb R^4)}
\le C\Big\|h_i\Big\|_{L^1(\mathbb R^4)}
$$
and
$$\Big\|II\Big\|_{L^{2,\infty}(\mathbb R^4)}
\le C\Big\||\Delta u_i|^2\Big\|_{L^1(\mathbb R^4)}
\le C\Big\|\nabla^2 u\Big\|_{L^2(\mathbb R^4)}^2.$$
For $III$, by integration by parts we have
\begin{eqnarray*}
&& \Big|III(x)\Big|\\
&=&\Big|\int_{\mathbb{R}^4}\nabla^2_y G(x-y)\Delta(|\nabla u_i|^2)u_i(y)\,dy\Big|\\
&\le & 2\left(\int_{\mathbb{R}^4}|\nabla^3_y G|(x-y)
|\nabla^2 u_i(y)||\nabla u_i(y)|\,dy
+\int_{\mathbb{R}^4}|\nabla^2_y
G|(x-y)|\nabla^2 u_i(y)||\nabla u_i(y)|^2\,dy\right)\\
&\le & CI_1(|\nabla^2 u_i||\nabla u_i|)(x)
+CI_2(|\nabla^2 u_i||\nabla u_i|^2)(x).
\end{eqnarray*}
Therefore, we have
\begin{eqnarray*}
\Big\|III\Big\|_{L^{2,\infty}(\mathbb R^4)}
&\le& C\left(\Big\||\nabla^2 u_i||\nabla u_i|\Big\|_{L^{\frac 43}(\mathbb{R}^4)}
+\Big\||\nabla^2 u_i||\nabla u_i|^2\Big\|_{L^1(\mathbb{R}^4)}\right)\\
&\le & C \Big\|\nabla^2 u_i\Big\|_{L^2(\mathbb{R}^4)}
\left(\Big\|\nabla u_i\Big\|_{L^{4}(\mathbb{R}^4)}
+\Big\|\nabla u_i\Big\|^2_{L^4(\mathbb{R}^4)}\right).
\end{eqnarray*}
Similar to $III$, we can estimate $IV$ by
$$
\Big\|IV\Big\|_{L^{2,\infty}(\mathbb R^4)}
\le C\Big\|\nabla^2 u_i\Big\|_{L^2(\mathbb{R}^4)}
\left(\Big\|\nabla^2 u_i\Big\|_{L^2(\mathbb R^4)}+\Big\|\nabla u_i\Big\|_{L^{4}(\mathbb{R}^4)}
+\Big\|\nabla u_i\Big\|^2_{L^4(\mathbb{R}^4)}\right).
$$
Putting the estimates of $I,II, III, IV$ together, we have
\begin{eqnarray}\label{vi-estimate3}
\Big\|\nabla^2 v_i\Big\|_{L^{2,\infty}(\mathbb{R}^4)}
&\le &
 C\Big[\Big\|h_i\Big\|_{L^1(\mathbb{R}^4)}+\Big\|\nabla^2 u_i\Big\|_{L^2(\mathbb R^4)}^2
\nonumber\\
&&+\Big\|\nabla^2 u_i\Big\|_{L^2(\mathbb{R}^4)}[\Big\|\nabla u_i\Big\|_{L^4(\mathbb{R}^4)}+
\Big\|\nabla u_i\Big\|_{L^4(\mathbb{R}^4)}^2]\Big]\nonumber\\
&\le&
C\Big[
\Big\|h\Big\|_{L^1(\mathcal B_i)}+\Big\|\nabla^2 u\Big\|_{L^2(\mathcal B_i)}^{\frac32}
+\Big\|\nabla^2 u\Big\|_{L^2(\mathcal B_i)}^2
\Big].
\end{eqnarray}
This clearly yields (\ref{vi-estimate}), since $\|\nabla^2 u\|_{L^2(\mathcal B_i)}\le\epsilon\le 1$.

Since $w_i$ is a biharmonic function, the standard  estimate and (\ref{vi-estimate3}) imply
that  $w_i\in C^\infty(\mathcal B_i)$, and
\begin{eqnarray}\label{wi-estimate3}
\sup_{x\in \mathcal A_i}\Big\{|x|^2\Big|\nabla^2 w_i\Big|(x)\Big\}
&\le& C
\Big\|\nabla^2  w_i\Big\|_{L^{2,\infty}(\mathcal B_i)}\le C\Big[\Big\|\nabla^2 u\Big\|_{L^2(\mathcal B_i)}+\Big\|\nabla^2 v_i\Big\|_{L^{2,\infty}(\mathcal B_i)}\Big]
\nonumber\\
&\le & C\Big[ \Big\|h\Big\|_{L^1(\mathcal B_i)} +\Big\|\nabla^2 u\Big\|_{L^2(\mathcal B_i)}\Big].
\end{eqnarray}
This yields (\ref{wi-estimate}).

It follows from the above Claim that
\begin{eqnarray*}
&&\sup_{\lambda\ge 1}\ \lambda^2\sum_{i=0}^{K-1}\Big|\Big\{x\in
B_{2^{i+1}\sqrt{\frac{\epsilon}{\lambda}}}\backslash B_{2^i\sqrt{\frac{\epsilon}{\lambda}}}
: |\nabla^2 u(x)|>\lambda\Big\}\Big|\\
&\le& \sup_{\lambda\ge 1}\ \lambda^2\sum_{i=0}^{K-1}\Big|\Big\{x\in\mathcal A_i
: |\nabla^2 v_i(x)|>\frac{\lambda}2\Big\}\Big|
+\sup_{\lambda\ge 1}\ \lambda^2\sum_{i=0}^{K-1}\Big|\Big\{x\in\mathcal A_i
: |\nabla^2 w_i(x)|>\frac{\lambda}2\Big\}\Big|\\
&\le& C\sum_{i=0}^{K-1}\Big\|\nabla^2 v_i\Big\|_{L^{2,\infty}(\mathcal A_i)}^2
+C\sup_{\lambda\ge 1}\lambda^{-2} \sum_{i=0}^{K-1}  \int_{\mathcal A_i}|\nabla^2 w_i|^4\,dx\\
&\le& C\sum_{i=0}^{K-1}\Big\|\nabla^2 v_i\Big\|_{L^{2,\infty}(\mathcal A_i)}^2
+C\sup_{\lambda\ge 1}\lambda^{-2} \sum_{i=0}^{K-1} \Big|\mathcal A_i\Big|\Big\|\nabla^2 w_i\Big\|_{L^\infty(\mathcal A_i)}^4\\
&\le& C \sum_{i=0}^{K-1}\left(\|h\|_{L^1(\mathcal B_i)}^2+\|\nabla^2 u\|_{L^2(\mathcal B_i)}^3\right)\\
&+&C\Big[\sum_{i=0}^{K-1}(2^i\sqrt{\epsilon})^{-4}\Big]
\cdot\sup_{0\le i\le K-1} [\|h\|_{L^1(\mathcal B_i)}^4+\|\nabla^2 u\|_{L^2(\mathcal B_i)}^4]\\
&=&I+II.
\end{eqnarray*}
Since H\"older inequality implies
$$\Big\|h\Big\|_{L^1(\mathcal B_i)}\le C\delta^{4(1-\frac1p)}
\Big\|h\Big\|_{L^p(B_{2\delta})},$$
we have
$$\sup_{0\le i\le K-1}\Big[ \|h\|_{L^1(\mathcal B_i)}^4+\|\nabla^2 u\|_{L^2(\mathcal B_i)}^4\Big]
\le \epsilon^4+C\Big\|h\Big\|_{L^p(B_{2\delta})}^4\delta^{16(1-\frac1p)}.$$
It is easy to see
\begin{eqnarray*}
I&\le& C\Big[\sup_{0\le i\le K-1}\|h||_{L^1(\mathcal B_i)}\Big]
\int_{\cup_{i=0}^{K-1}\mathcal B_i}|h|
+C\epsilon \int_{\cup_{i=0}^{K-1}\mathcal B_i}|\nabla^2 u|^2\\
&\le& C\delta^{4(1-\frac1{p})}\|h\|_{L^p(B_{2\delta})} \|h\|_{L^1(B_{2\delta})}
+C\epsilon\|\nabla^2 u\|_{L^2(B_{2\delta})}^2\\
&\le& C\delta^{8(1-\frac1{p})}\|h\|_{L^p(B_{2\delta})}^2
+C\epsilon\|\nabla^2 u\|_{L^2(B_{2\delta})}^2,
\end{eqnarray*}
and
\begin{eqnarray*}
II&\le& C\epsilon^{-2} \sup_{0\le i\le K-1} \Big[\|h\|_{L^1(\mathcal B_i)}^4+\|\nabla^2 u\|_{L^2(\mathcal B_i)}^4\Big]\\
&\le& C\Big(\epsilon^2+\epsilon^{-2}\delta^{16(1-\frac1p)}\Big\|h\Big\|_{L^p(B_{2\delta})}^4\Big).
\end{eqnarray*}
Putting these estimates together we can obtain
\begin{eqnarray*}
\Big\|\nabla^2 u\Big\|_{L^{2,\infty}(B_\delta\setminus B_{2Rr})}^2
&\le& C\Big(\epsilon^2+\epsilon\|\nabla^2 u\|_{L^2(B_{2\delta})}^2+\delta^4
+\delta^{8(1-\frac1{p})}\|h\|_{L^p(B_{2\delta})}^2\\
&&+
\epsilon^{-2}\delta^{16(1-\frac1p)}\Big\|h\Big\|_{L^p(B_{2\delta})}^4\Big).
\end{eqnarray*}
This clearly implies (\ref{small_2infty}). The proof is complete. \qed\\

\section{Proof of Theorem \ref{energy_identity} and \ref{energy_identity1}}
\setcounter{equation}{0}
\setcounter{theorem}{0}

This section is devoted to the proof of ``$=$'' of (\ref{energy_id1}). The argument is based on the duality between
$L^{2,1}$ and $L^{2,\infty}$.

\medskip
\noindent{\bf Competition of Proof of Theorem} \ref{energy_identity}:

\smallskip
For simplicity, we may assume $\Sigma=\{0\}\subset\Omega$ is a single point. In particular, $u_m\rightarrow u$ in
$W^{2,2}_{\rm{loc}}(B_{r_1}\setminus\{0\})$ for some $r_1>0$.
By an induction argument similar to that of \cite{DT} in the context of harmonic maps, we may assume that there is
only one bubble in $ B_{r_1}$, i.e. $L_1=1$. Then for any
$\epsilon>0$, there exist $r_m\downarrow 0$, $R\ge 1$ sufficiently large, and $0<\delta\le \epsilon^{\frac{p}{4(p-1)}}$
such that  for $m$ sufficiently large, it holds
\begin{equation}\label{dy_small}
\int_{B_{2\rho}\setminus B_\rho} |\nabla^2 u_m|^2
\le \epsilon^2, \ \forall Rr_m\le\rho\le 2\delta.
\end{equation}
Therefore by Lemma 4.1 we have
\begin{eqnarray}\label{2infty_norm}
\Big\|\nabla^2 u_m\Big\|_{L^{2,\infty}(B_\delta\setminus B_{2Rr_m})}
&\le& C\Big[\epsilon+\delta^2+\delta^{4(1-\frac{1}{p})}
\Big\|h_m\Big\|_{L^p(B_{2\delta})}+\epsilon^\frac12
\Big\|\nabla^2 u_m\Big\|_{L^2(B_{2\delta})}\nonumber\\
&&+\epsilon\Big\|h_m\Big\|_{L^p(B_{2\delta})}^2\Big]\nonumber\\
&\le & C\epsilon^\frac12.
\end{eqnarray}

On the other hand, by Theorem 2.3 we have
\begin{equation}\label{L21_bound}
\Big\|\nabla^2 u_m\Big\|_{L^{2,1}(B_\delta)}
\le C\Big(\Big\|\nabla^2 u_m\Big\|_{L^2(B_{2\delta})}
+\delta^{4(1-\frac1{p})}\Big\|h_m\Big\|_{L^p(B_{2\delta})}\Big)
\le C.
\end{equation}
Therefore by the duality between $L^{2,1}$ and $L^{2,\infty}$,
we have
\begin{eqnarray}
\Big\|\nabla^2 u_m\Big\|^2_{L^{2}(B_{\delta}\backslash B_{2Rr_m})} &\le &
C\Big\|\nabla^2 u_m\Big\|_{L^{2,1}(B_{\delta}\backslash B_{2Rr_m})}
\Big\|\nabla^2 u_m\Big\|_{L^{2,\infty}(B_{\delta}\backslash B_{2Rr_m})}\nonumber\\
 &\le &
C\epsilon^{\frac12}.\label{L2_bound}
\end{eqnarray}
Since $\epsilon>0$ is arbitrary, this yields that (\ref{energy_id1}) holds. It is easy to see
that (\ref{energy_id2}) follows from (\ref{energy_id1})
and the pointwise inequality $|\nabla u_m|^2\le |\Delta u_m|$.
The proof of Theorem \ref{energy_identity} is now complete.
\qed

\begin{remark} Under the same assumption as Theorem \ref{energy_identity} and notations as above, it holds for any $x_0\in\Sigma$,
\begin{equation}\label{no_oscillation}
\lim_{\delta\downarrow 0}\lim_{R\uparrow \infty}\lim_{m\rightarrow
\infty} {osc}_{B_{\frac{\delta}2}(x_0)\setminus B_{2Rr_m}(x_0)} u_m=0.
\end{equation}
\end{remark}
\pf By (\ref{L2_bound}), we have that for any
$\epsilon>0$ there are $R>1$ sufficiently large, $\delta>0$
sufficiently small, and $m\ge 1$ sufficiently large so that
(\ref{L2_bound}) holds.
By Fubini's theorem, we may assume
$$Rr_m\int_{\partial B_{2Rr_m}(x_0)}|\nabla^2 u_m|^2\le C\epsilon^\frac12.$$
Let $v_m: B_{2Rr_m}\to\mathbb S^k$ be a minimizing biharmonic map
extension of $u_m$ such that $(v_m,\nabla v_m)=(u_m, \nabla u_m) \ {\rm{on}}\ \partial B_{2Rr_m}(x_0).$
Then we would have (cf. \cite{scheven2})
$$\int_{B_{2Rr_m}(x_0)}|\nabla^2 v_m|^2\le
C Rr_m\int_{\partial B_{2Rr_m}(x_0)}|\nabla^2 u_m|^2\le C\epsilon^\frac12.$$
Now we define $w_m: B_\delta(x_0)\to\mathbb S^{k}$ by
letting
$$w_m=\begin{cases} v_m & {\rm{in}}\ B_{2Rr_m}(x_0)\\
u_m & {\rm{in}}\ B_\delta(x_0)\setminus B_{2Rr_m}(x_0).
\end{cases}
$$
Then $w_m$ is an approximate biharmonic map to $\mathbb S^k$ with
bi-tension field $\widetilde {h_m}$ given by
$$\widetilde{h_m}=\begin{cases}0 & \ {\rm{in}}\ B_{2Rr_m}(x_0)\\
 h_m & \ {\rm{in}}\ B_{\delta}(x_0)\setminus B_{2Rr_m}(x_0).
\end{cases}
$$
Therefore by Theorem 2.3 we have
$$
\Big\|\nabla^2 w_m\Big\|_{L^{2,1}(B_\delta(x_0))}^2\le C\epsilon^\frac12.$$
Since there exists a harmonic function $\widetilde {w_m}$ on $B_\delta(x_0)$, with $\|\nabla^2 \widetilde {w_m}\|_{L^2(B_\delta(x_0))}\le
C\epsilon^\frac12$,
such that
$$w_m(x)=\int_{B_\delta(x_0)}|x-y|^{-2}\Delta w_m(y)\,dy
+\widetilde {w_m}(x), \ x\in B_{\delta}(x_0).$$
It is easy to see
$$\hbox{osc}_{B_{\frac{\delta}2}(x_0)} \widetilde{w_m}
\le C\epsilon^\frac12,$$
while by the duality between $L^{2,1}$ and $L^{2,\infty}$ we have
\begin{eqnarray*}
&&\Big|\int_{B_{\delta}(x_0)}|x-y|^{-2}
\Delta w_m(y)\,dy\Big|_{L^\infty(B_{\frac{\delta}2}(x_0))}
\\&\le & C
\Big\|\frac1{|x|^2}\Big\|_{L^{2,\infty}(B_1)}
\Big\|\nabla^2 w_m\Big\|_{L^{2,1}(B_\delta(x_0))}\le C\epsilon^\frac14.
\end{eqnarray*}
Thus (\ref{no_oscillation}) follows.
\qed\\

\noindent{\bf Proof of Theorem \ref{energy_identity1}}:

\smallskip
It follows from the energy inequality (1.8)
that there exists $t_m\uparrow +\infty$ such that
$u_m(\cdot)=u(\cdot,t_m)$ is an approximate biharmonic map into
$\mathbb S^k$ with bi-tension fields
$h_m=u_t(\cdot, t_m)\in L^2(\Omega)$ satisfying
$$\Big\|h_m\Big\|_{L^2(\Omega)}
=\Big\|u_t(\cdot,t_m)\Big\|_{L^2(\Omega)}\rightarrow 0.$$
Moreover,
$$\Big\|u_m\Big\|_{W^{2,2}(\Omega)}\le
C\Big\| u_0\Big\|_{W^{2,2}(\Omega)}.$$
Therefore we may assume that after taking another subsequence,
$u_m\rightharpoonup u_\infty$ in $W^{2,2}(\Omega,\mathbb S^k)$.
It is easy to see that $u_\infty$ is a biharmonic map so that
$u_\infty\in C^\infty(\Omega,\mathbb S^k)$ (see \cite{W1}).
All other conclusions follow directly from Theorem \ref{energy_identity}. \qed

\end{document}